\ifdefined\journalstyle

 \documentclass[graybox]{svmult}

 \usepackage{type1cm}        
 \usepackage{makeidx}        
 \usepackage{graphicx}       
 \usepackage{multicol}        
 \usepackage[bottom]{footmisc}

 \usepackage{newtxtext}       %
 \usepackage[varvw]{newtxmath}       

 \input{common-macros}

\else
 \documentclass[11pt,a4paper]{article}
 \usepackage{latexsym,amsfonts,amsmath,amsthm,bm,graphics,multirow,microtype}

 \newtheorem{theorem}{Theorem}
 \newtheorem{lemma}[theorem]{Lemma}
 \newtheorem{corollary}[theorem]{Corollary}
 \usepackage[usenames,dvipsnames,svgnames,table]{xcolor}

 \usepackage[width=.75\textwidth]{caption}

\usepackage{graphicx}       

\usepackage{amsmath,amsfonts,amssymb,bm,amsthm,mathtools,mathrsfs}

\def\citep#1#2{\cite[{#1}]{#2}}

\usepackage{algorithm}
\usepackage{algorithmic}

\newcommand{\bbN}{{\mathbb{N}}}

\newcommand{\bbR}{{\mathbb{R}}}

\newcommand{\N}{{\mathbb{N}}} 

\newcommand{\ind}{{\rm ind}}    

\newcommand{\bsa}{{\boldsymbol{a}}}

\newcommand{\bsf}{{\boldsymbol{f}}}

\newcommand{\bsm}{{\boldsymbol{m}}}

\newcommand{\bst}{{\boldsymbol{t}}}

\newcommand{\bsx}{{\boldsymbol{x}}}
\newcommand{\bsy}{{\boldsymbol{y}}}
\newcommand{\bsz}{{\boldsymbol{z}}}

\newcommand{\bsY}{{\boldsymbol{Y}}}
\newcommand{\bsZ}{{\boldsymbol{Z}}}

\newcommand{\bszero}{{\boldsymbol{0}}} 

\newcommand{\bsnu}{{\boldsymbol{\nu}}}

\newcommand{\calK}{{\mathcal{K}}}

\newcommand{\calO}{{\mathcal{O}}}

\newcommand{\setu}{{\mathfrak{u}}}

\newcommand{\rd}{\,\mathrm{d}} 

\fi

\usepackage{pdfsync} 

\definecolor{darkgreen}{RGB}{0,100,0}
\definecolor{darkmagenta}{RGB}{139,0,139}

\usepackage{scrextend}

\begin{document}

\ifdefined\journalstyle

 \title*{Lattice-based kernel approximation and serendipitous weights for parametric PDEs in very high dimensions}

 \titlerunning{Lattice-based kernel approximation and serendipitous weights in very high dimensions}

 \author{Vesa Kaarnioja\and Frances Y. Kuo\and Ian H. Sloan}

 \institute{V. Kaarnioja \at Department of Mathematics and Computer Science, Free
  University of Berlin, Arnimallee 6, 14195 Berlin, Germany, \email{vesa.kaarnioja@fu-berlin.de}
  \and
  F. Y. Kuo \at School of Mathematics and Statistics, UNSW Sydney, Sydney NSW 2052, Australia,
  \email{f.kuo@unsw.edu.au}
  \and
  I. H. Sloan \at School of Mathematics and Statistics, UNSW Sydney, Sydney NSW 2052, Australia,
  \email{i.sloan@unsw.edu.au}}

\else

 \title{Lattice-based kernel approximation and \\
 serendipitous weights for parametric PDEs in \\ very high dimensions}

 \author{Vesa Kaarnioja\footnote{Department of Mathematics and Computer Science, Free
  University of Berlin, Arnimallee 6, 14195 Berlin, Germany. \newline Email: vesa.kaarnioja@fu-berlin.de}
  \and Frances Y. Kuo\footnote{School of Mathematics and Statistics, UNSW Sydney, Sydney NSW 2052,
 Australia. Email: f.kuo@unsw.edu.au}
  \and Ian H. Sloan\footnote{School of Mathematics and Statistics, UNSW Sydney, Sydney NSW 2052,
 Australia. Email: i.sloan@unsw.edu.au}
 }
 \date{October 2026}

\fi

\maketitle

\abstract{We describe a fast method for solving elliptic partial
differential equations (PDEs) with uncertain coefficients using kernel
interpolation at a lattice point set. By representing the input random
field of the system using the model proposed by Kaarnioja, Kuo, and Sloan
(SIAM J.~Numer.~Anal.~2020), in which a countable number of independent
random variables enter the random field as periodic functions, it was
shown by Kaarnioja, Kazashi, Kuo, Nobile, and Sloan (Numer.~Math.~2022)
that the lattice-based kernel interpolant can be constructed for the PDE
solution as a function of the stochastic variables in a highly efficient
manner using fast Fourier transform (FFT). In this work, we discuss the
connection between our model and the popular ``affine and uniform model''
studied widely in the literature of uncertainty quantification for PDEs
with uncertain coefficients. We also propose a new class of product
weights entering the construction of the kernel interpolant, which
dramatically improve the computational performance of the kernel
interpolant for PDE problems with uncertain coefficients, and allow us to
tackle function approximation problems up to very high dimensionalities.
Numerical experiments are presented to showcase the performance of the new
weights.}

\section{Introduction}

As in the major survey~\cite{CD15}, we consider parametric elliptic
partial differential equations (PDEs) of the form
\begin{equation}\label{eq:pde}
 -\nabla\cdot (\widetilde{a}(\bsx,\bsz)\nabla \widetilde{u}(\bsx, \bsz)) = q(\bsx),
 \quad \bsx \in D \subset
 \mathbb{R}^d, \quad \bsz \in [-1,1]^s,
\end{equation}
subject to the homogeneous Dirichlet boundary condition
$\widetilde{u}(\bsx,\bsz) = 0$ for $\bsx \in \partial D$, where
$D\subset\mathbb R^d$ is a bounded Lipschitz domain, with $d$ typically
$2$ or $3$, and $\widetilde{a}$ is an uncertain input field given in terms
of parameters $\bsz = (z_1,z_2,\ldots, z_s)\in [-1,1]^s$ by
\begin{equation}\label{eq:affine}
 \widetilde{a}(\bsx, \bsz) := a_0(\bsx) +  \sum_{j=1}^s z_j\, \psi_j(\bsx), \quad \bsx \in D,
 \quad \bsz \in [-1,1]^s.
\end{equation}
Here the parameters $z_j$ represent independent random variables
distributed on the interval $[-1,1]$ according to a given probability
distribution with density $\rho$.  The essential feature giving rise to
major difficulties is that $s$ may be very large.

We assume that the univariate density $\rho$ has the Chebyshev or arcsine
form
\begin{equation}\label{eq:arcsin}
 \rho(z):=\frac{1}{\pi \sqrt{1-z^2}},
\quad z \in [-1,1].
\end{equation}
This is one of only two specific probability densities considered in the
recent monograph~\cite{ABS22}, the other being the constant density
$\rho(z) = \frac{1}{2}$. In the method of generalized polynomial chaos
(GPC), see~\cite{XK02}, the density \eqref{eq:arcsin} is associated with
Chebyshev polynomials of the first kind as the univariate basis functions.
It might be argued that in many applications the choice between these two
densities is a matter of taste rather than conviction.

In Section~\ref{sec:transform} we set $\bsz = \sin(2\pi \bsy)$
component-wise, and so recast the problem
\eqref{eq:pde}--\eqref{eq:affine} with the density~\eqref{eq:arcsin} to
one in which the dependence on a new stochastic variable $\bsy$ is
periodic, and the solution becomes
\begin{equation}\label{eq:u}
u(\bsx, \bsy):= \widetilde{u}(\bsx,\sin(2\pi\bsy)), \quad \bsx \in D, \quad\bsy \in [0,1]^s\, =:\, U_s.
\end{equation}
More precisely, we show that $u(\bsx, \bsy)$ satisfies
\begin{equation}\label{eq:notilde_pde}
 -\nabla\cdot (a(\bsx,\bsy)\nabla u(\bsx, \bsy)) = q(\bsx), \quad \bsx \in D \subset
\mathbb{R}^d, \quad \bsy \in U_s
\end{equation}
 and that
\eqref{eq:affine} together with the probability density \eqref{eq:arcsin}
can be replaced by
\begin{equation}\label{eq:periodic}
 a(\bsx, \bsy) := \widetilde{a}(\bsx,\sin(2\pi\bsy))
 = a_0(\bsx) +  \sum_{j=1}^s  \sin(2\pi y_j)\, \psi_j(\bsx), \quad \bsx \in D, \quad\bsy \in U_s,
\end{equation}
with $\bsy = (y_1, y_2,\ldots, y_s)$, where the parameters $y_j$ represent
i.i.d.~random variables \emph{uniformly} distributed on $[0,1]$. The
equivalence of \eqref{eq:affine} subject to density \eqref{eq:arcsin} on
the one hand, with \eqref{eq:periodic} subject to uniform density on the
other, is meant in the sense that both fields have exactly the same joint
probability distribution in the  domain $D$; for a precise statement see
Theorem~\ref{theorem:22} and Corollary~\ref{corollary:23} below.

In Section~\ref{section:kernel} we describe the lattice-based kernel
interpolant of~\cite{KKKNS22} and discuss its properties, error bounds and
cost of construction and evaluation. The essence of the method is that the
dependence on the parameter $\bsy$ is approximated by a linear combination
of periodic kernels, each with one leg fixed at a point $\bst_k, k =
1,\ldots,n$, of a carefully chosen $s$-dimensional rank-$1$ lattice, with
coefficients fixed by interpolation. The kernel is the reproducing kernel
of a weighted $s$-variate Hilbert space $H$ of dominating mixed smoothness
$\alpha\in\bbN$. (The parameter labelled $\alpha$ in \cite{KKKNS22} is
here replaced by $2\alpha$, to make $\alpha$ correspond to the order of
mixed derivatives, see the norm~\eqref{eq:norm} below.)

In Section~\ref{sec:pde} we summarize the results from the
paper~\cite{KKKNS22} on applying the lattice-based kernel interpolant to
the PDE problem \eqref{eq:notilde_pde}--\eqref{eq:periodic}. The main
result is that, provided the fluctuation coefficients $\psi_j$ in
\eqref{eq:periodic} decay suitably, the $L_2$ error (with respect to both
$\bsx$ and $\bsy$) of the kernel interpolant to the PDE solution is of the
order
\begin{equation}\label{eq:errorbd}
 \calO(n^{-\alpha/2+\delta}), \quad \delta >0,
\end{equation}
where the implied constant depends on $\delta$ but is independent of ~$s$.
The convergence rate $n^{-\alpha/2}$ is known (see \cite{BKUV17}) to be
the best that can be achieved in the sense of worst case $L_2$ error for
any approximation that (as in \cite{KKKNS22} and here) uses information
only at the points of a rank-$1$ lattice.

The present paper improves upon the method as presented in \cite{KKKNS22}
in two different ways: first by making the method much more efficient; and
second by making it much more accurate in difficult cases.  The combined
effect is to greatly increase the range of dimension $s$ that are
realistically achievable. We demonstrate this by carrying out explicit
nontrivial calculations with $s= 1000$, compared with $s = 100$ in
\cite{KKKNS22}.

Both benefits are achieved in the present paper by the introduction of a
new kind of ``weights''. As mentioned earlier, our approximation makes
explicit use of the reproducing kernel of a certain weighted Hilbert space
$H$. The definition of that Hilbert space involves positive numbers
$\gamma_\setu$ called weights, see \eqref{eq:norm} below. The weights
naturally occur also in the kernel, see \eqref{eq:ker}. The best
performing weights found in~\cite{KKKNS22} were of the so-called ``SPOD''
form (see \eqref{eq:spod} below). For SPOD weights the cost of evaluating
the kernel interpolant is high, because the cost of kernel evaluation is
high.

In Section~\ref{sec:serendipity} we introduce novel weights, which are
``product'' weights rather than SPOD weights, allowing the kernels to be
evaluated at a cost of order $s$. We still have a rigorous error bound of
the form \eqref{eq:errorbd}. The theoretical downside is that the implied
constant is no longer independent of $s$. We refer to these weights as
``serendipitous'', with the word ``serendipity'' meaning ``happy discovery
by accident''.

In Section~\ref{sec:numex} we give numerical results, which show that with
the lattice-based kernel interpolant from~\cite{KKKNS22} and these
weights, problems with dimension $s$ of a thousand or more can easily be
studied. For not only are they cheaper and easier to use than SPOD
weights, but also in difficult problems they lead empirically (and to our
surprise) to much smaller errors, while producing similar errors to SPOD
weights for easier problems.

A different kind of product weight was developed in \cite{KKKNS22} by
adhering to the requirement that the error bound be independent of
dimension, but those weights were found to have limited applicability, and
did not show the remarkable performance reported here.

Many other methods have been proposed for $L_2$ approximation in the
multivariate setting. The review article~\cite{CD15} and the
monograph~\cite{ABS22} take the approximating function to be a
multivariate polynomial; as a result a major part of their analysis is
inevitably concerned with sparsification of the basis set, since otherwise
the ``curse of dimensionality'' would preclude large values of $s$.
Recently other methods have been proposed \cite{BKPU22,BKUV17,KPV15,KV19},
some of which are optimal with respect to order of convergence, in the
sense of producing error bounds of order
\begin{align} \label{eq:order}
 \calO(n^{-\alpha}) \quad \mbox{or} \quad  \calO(n^{-\alpha} (\log n)^{\beta})
\end{align}
for some $\beta$, a rate with respect to the exponent of $n$ that is the
same as for approximation numbers, see \cite{KU21,NSU22,DKU23}. Obviously
\eqref{eq:order} displays a better convergence rate than
\eqref{eq:errorbd}, but it has yet to be demonstrated that any of these
methods has the potential for solving in practice the very
high-dimensional problems seen in the numerical calculations of this
paper.

The application of lattice point sets together with kernel interpolation
has gained a lot of attention in the recent years. The paper~\cite{ZLH06}
appears to have been the first to consider this approach, later the
paper~\cite{ZKH09} obtained dimension-independent error estimates in
weighted spaces using product weights.

The use of lattice points for approximation has been facilitated by the
development of fast component-by-component (CBC) constructions for
lattices under different assumptions on the form of the weights,
see~\cite{KSW06,CKNS20,CKNS21}. This enables the generation of tailored
lattices for large-scale high-dimensional approximation problems such as
those arising in uncertainty quantification.

This paper is organized as follows. In Section~\ref{sec:transform}, we
describe the connection between the so-called ``affine model'' and the
``periodic model'' introduced in~\cite{KKS20}. The lattice-based kernel
interpolant of~\cite{KKKNS22} is summarized in
Section~\ref{section:kernel}. The application to parametric elliptic PDEs
with uncertain coefficients is discussed in Section~\ref{sec:pde}. The new
class of serendipitous weights for the construction of the kernel
interpolant for parametric PDE problems is introduced in
Section~\ref{sec:serendipity}. Numerical experiments assessing the
performance of those weights are presented in Section~\ref{sec:numex}. The
paper ends with some conclusions.

\section{Transforming to the periodic setting}\label{sec:transform}

The equivalence of the affine probability model given by \eqref{eq:affine}
and \eqref{eq:arcsin} with the periodic formulation in \eqref{eq:periodic}
is expressed in Corollary~\ref{corollary:23} below. It states that the
probability distributions in the two cases are identical. A more general
result is stated in Theorem~\ref{theorem:22}. It rests in turn on the
following lemma, stating that if $Y$ is a random variable uniformly
distributed on $[0,1]$, then $\sin(2\pi Y)$ has exactly the same
probability distribution as a random variable distributed on $[-1,1]$ with
the arcsine probability density~\eqref{eq:arcsin}. While the proof is
elementary, it has one slightly unusual feature, that the change of
variable from $y \in [0,1]$ to $z = \sin(2\pi y)$ is not monotone.

\begin{lemma}\label{lemma:21}
Let $Z$ be a random variable distributed on $[-1,1]$ with density $\rho$
given by \eqref{eq:arcsin}, and let $Y$ be a random variable uniformly
distributed on $[0,1]$. Then for all $t\in \mathbb{R}$ we have
\[
 \mathbb{P}[Z\le t] = \mathbb{P}[\sin(2\pi Y)\le t].
\]
\end{lemma}

\begin{proof}
 We first write the left-hand side explicitly as an integral:
 \[
 \mathbb{P}[Z\le t] = \int_{-1}^1 \ind(t-z)\,\frac{1}{\pi \sqrt{1 - z^2}}\rd z,
 \]
where $\ind$ is the indicator function, taking the value $1$ if the
argument is non-negative, and the value $0$ if it is negative.  The next
step is to make a change of variable $z = \sin(2\pi y)$, but note that
this is not permissible for $y$ on the whole interval $[0,1]$ because
$\sin(2\pi y)$ is not monotone. Noting that $\sin(2\pi y)$ is
$1$-periodic, it is sufficient to consider $y$ in the two intervals
$[-1/4, 1/4]$ and $[1/4, 3/4]$ separately (the point being that $\sin(2\pi
y)$ is monotone in each subinterval, and that together the two
subintervals  make a full period). For the first subinterval we have
\[
 z = \sin(2\pi y), \quad y = \frac{\arcsin(z)}{2\pi} \in [-1/4,1/4],
 \]
so that we obtain
\[
 \mathbb{P}[Z\le t] = \int_{-1/4}^{1/4} \ind(t-\sin(2\pi y))
 \frac{2 \pi \cos(2\pi y)}{\pi \cos(2\pi y)}\rd y
  =\int_{-1/4}^{1/4} \ind(t-\sin(2\pi y))\,2\rd y.
 \]

Similarly, for the second interval $[1/4,3/4]$ we have
 \[
 z = \sin(2\pi y), \quad y = \frac{1}{2} - \frac{\arcsin(z)}{2\pi} \in [1/4,3/4],
 \]
 so that
 \[
 \mathbb{P}[Z\le t]  =\int_{1/4}^{3/4} \ind(t-\sin(2\pi y))\,2 \rd y.
 \]
Adding the two results together, and dividing by 2, we obtain
\[
 \mathbb{P}[Z\le t] =\int_{-1/4}^{3/4} \ind(t-\sin(2\pi y))\rd y
 =\int_{0}^{1} \ind(t-\sin(2\pi y)) \rd y
 = \mathbb{P}[\sin(2\pi Y)\le t],
\]
with the second equality following from the periodicity of $\sin(2\pi y)$.
This completes the proof of the lemma.
\end{proof}

\begin{theorem}\label{theorem:22}
Let $\widetilde{Q}(\bsZ) = \widetilde{Q}(Z_1, Z_2, \ldots,Z_s)$ be a
real-valued random variable that depends on $s$ i.i.d.~real-valued random
variables $Z_1, Z_2, \ldots, Z_s$, where each $Z_j$ is distributed on
$[-1,1]$ with density $\rho$ given by~\eqref{eq:arcsin}. Let
$Q(\bsY)=Q(Y_1,Y_2,\ldots, Y_s)$ be another random variable defined by
\[ 
 Q(\bsY) = \widetilde{Q}(\sin(2\pi\bsY)),
\] 
where the $Y_j$ are i.i.d.~uniformly distributed random variables on
$[0,1]$. Then for all $q\in \mathbb{R}$ we have
\[ 
 \mathbb{P}[Q(\cdot) \le q] =  \mathbb{P}[\widetilde{Q}(\cdot) \le q],
\] 
where the probability on the left is with respect to uniform density, while that on the right is with respect to a product of the univariate densities \eqref{eq:arcsin}.
\end{theorem}

\begin{proof}
This follows immediately by applying  Lemma~\ref{lemma:21} to each random
variable, together with Fubini's theorem.
\end{proof}

It follows from the theorem that the input parametric field given by
\eqref{eq:affine} subject to the density \eqref{eq:arcsin} can with equal
validity be expressed as \eqref{eq:periodic}, with each $y_j$ uniformly
distributed on $[0,1]$. We state this as a corollary:

\begin{corollary}\label{corollary:23} Let $\bsZ$ and $\bsY$ be as in
Theorem \ref{theorem:22}, and $\widetilde{u}(\bsx, z)$ and $u(\bsx,\bsy)$
be as in~\eqref{eq:pde} and~\eqref{eq:u}.  The random variable
$\widetilde{u}(\bsx,\bsZ)$ for $\bsx\in D$ has the same joint probability
distribution with respect to the product density $\prod_{j=1}^s \rho(z_j)$
as $u(\bsx, \bsY)$ has with respect to the uniform product density.
\end{corollary}

This is the periodic formulation introduced by~\cite{KKS20}, except for
the trivial difference of a different normalising factor: in~\cite{KKS20}
the sum was multiplied by $1/\sqrt{6}$ to ensure that the variance of the
field was the same as that of a uniformly distributed affine variable
defined on $[-1/2,1/2]$. In effect we are here redefining the fluctuation
coefficients $\psi_j$.

\section{The kernel interpolant}\label{section:kernel}

We assume that $f: [0,1]^s \to \bbR$ belongs to the \emph{weighted
periodic unanchored Sobolev space} $H$ of dominating mixed smoothness of
order $\alpha\in\bbN$, with norm
\begin{align}\label{eq:norm}
 \|f\|_H^2 :=
 \sum_{\setu\subseteq\{1:s\}}\frac{1}{(2\pi)^{2\alpha |\setu|} \gamma_\setu}
 \int_{[0,1]^{|\setu|}}\!\!\;\bigg|\int_{[0,1]^{s-|\setu|}} \bigg(\prod_{j\in\setu}
 \frac{\partial^{\alpha}}{\partial y_j^{\alpha}}\bigg)
 f(\bsy)\rd\bsy_{-\setu}\bigg|^2\,\rd\bsy_\setu,
\end{align}
where $\{1:s\} := \{1,2,\ldots,s\}$, $\bsy_{\setu}:=(y_j)_{j\in\setu}$ and
$\bsy_{-\setu}:=(y_j)_{j\in\{1:s\}\setminus\setu}$. The inner product
$\langle \cdot,\cdot\rangle_H$ is defined in an obvious way. Here we have
replaced the traditional notation of the smoothness parameter $\alpha$ by
$2\alpha$, so that our $\alpha$ corresponds exactly to the order of
derivatives in the norm \eqref{eq:norm}. The space $H$ is a special case
of the weighted Korobov space which has a real smoothness parameter
$\alpha$ characterizing the rate of decay of Fourier coefficients, see
e.g., \cite{SW01,KSW06,CKNS20,CKNS21}.

The parameters $\gamma_{\setu}$ for $\setu\subseteq\{1:s\}$ in the norm
\eqref{eq:norm} are \emph{weights} that are used to moderate the relative
importance between subsets of variables in the norm, with
$\gamma_{\emptyset}:=1$. There are $2^s$ weights in full generality, far
too many to prescribe one by one. In practice we must work with weights of
restricted forms. The following forms of weights have been considered in
the literature:
\begin{itemize}
\item \emph{Product weights} \cite{SW01}: $\gamma_\setu =
    \prod_{j\in\setu} \gamma_j$, specified by one sequence
    $(\gamma_j)_{j\ge 1}$.
\item \emph{POD weights} (product and order dependent) \cite{KSS12}:
    $\gamma_\setu = \Gamma_{|\setu|} \prod_{j\in\setu} \gamma_j$,
    specified by two sequences $(\Gamma_\ell)_{\ell\ge 0}$ and
    $(\gamma_j)_{j\ge 1}$.
\item \emph{SPOD weights} (smoothness-driven product and order
    dependent) \cite{DKLNS14}: $\gamma_\setu =
    \sum_{\bsnu_\setu\in\{1:\alpha\}^{|\setu|}} \Gamma_{|\bsnu_\setu|}
    \prod_{j\in\setu} \gamma_{j,\nu_j}$, specified by the sequences
    $(\Gamma_\ell)_{\ell\ge 0}$ and $(\gamma_{j,\nu})_{j\ge 1}$ for
    each $\nu=1,\ldots,\alpha$, where $|\bsnu_\setu| :=
    \sum_{j\in\setu} \nu_j$.
\end{itemize}

The space $H$ is a reproducing kernel Hilbert space (RKHS) with the kernel
\begin{equation}\label{eq:ker}
K(\bsy,\bsy'):=\sum_{\setu\subseteq\{1:s\}}\gamma_{\setu}\prod_{j\in\setu}
\eta_{\alpha}(y_j,y_j'),\quad \bsy,\bsy'\in [0,1]^s,
\end{equation}
where
$$
\eta_{\alpha}(y,y'):=\frac{(2\pi)^{2\alpha}}{(-1)^{\alpha+1}(2\alpha)!}B_{2\alpha}(\{y-y'\}),\quad y,y'\in[0,1],
$$
with $B_{\ell}(y)$ denoting the Bernoulli polynomial of degree $\ell$, and
the braces $\{\cdot\}$ denoting the fractional part of each component of
the argument. As concrete examples, we have $B_2(y) = y^2-y + 1/6$ and
$B_4(y) = y^4-2y^3+y^2-1/30$. The kernel $K$ is easily seen to satisfy the
two defining properties of a reproducing kernel, namely that
$K(\cdot,\bsy)\in H$ for all $\bsy\in [0,1]^s$ and $\langle f,
K(\cdot,\bsy)\rangle_H = f(\bsy)$ for all $f\in H$ and all $\bsy\in
[0,1]^s$.

For a given $f\in H$, we consider an approximation of the form
\begin{equation} \label{eq:appro}
 A_n^*(f)(\bsy):=f_n(\bsy)
 :=\sum_{k=1}^n a_k\, K(\bst_k,\bsy),\quad \bsy\in [0,1]^s,
\end{equation}
where $a_1, \ldots, a_n\in\mathbb{R}$ and the nodes
\begin{align} \label{eq:latticepoints}
\bst_k:=\bigg\{\frac{k\bsz_{\rm gen}}{n}\bigg\}\quad\text{for}~k = 1,\ldots,n
\end{align}
are the points of a rank-1 lattice, with $\bsz_{\rm
gen}\in\{1,\ldots,n-1\}^s$ being the \emph{lattice generating vector}.
Since the kernel $K$ is periodic, the braces $\{\cdot\}$ indicating the
fractional part in the definition of $\bst_k$ can be omitted when we
evaluate the kernel.

The \emph{kernel interpolant} $f_n\in H$ is a function of the
form~\eqref{eq:appro} that interpolates $f$ at the lattice points,
\[ 
 f_n(\bst_{k'})=f(\bst_{k'})\quad\text{for all}~k'=1,\ldots,n\label{eq:interp}
\] 
with the coefficients $a_k$ therefore satisfying the linear system
\begin{align} \label{eq:linearsys}
\sum_{k=1}^n\calK_{k,k'}\,a_k=f(\bst_{k'})\quad\text{for all}~k'=1,\ldots,n,
\end{align}
where
\begin{align} \label{eq:matrix}
  \calK_{k,k'} = \calK_{k',k} := K(\bst_k,\bst_{k'})
  = K\left(\frac{(k-k')\bsz_{\rm gen}}{n},\mathbf 0\right)
 \quad\mbox{for } k,k'=1,\ldots,n.
\end{align}
The solution of the square linear sytem \eqref{eq:linearsys} exists and is
unique because $K$, by virtue of being a reproducing kernel, is positive
definite.

Moreover, because the nodes form a lattice, the matrix $\calK$ is a
\emph{circulant}, thus the coefficient vector $\bsa:=[a_1,\ldots,a_n]^{\rm
T}$ in~\eqref{eq:linearsys} can be solved in $\mathcal O(n\log n)$ time
using the fast Fourier transform (FFT) via
$$
\bsa={\tt ifft}({\tt fft}(\bsf)\,./\,{\tt fft}(\calK_{:,1})),
$$
where\, $./$\, indicates component-wise division of two vectors, $\bsf :=
[f(\bst_1),\ldots,f(\bst_n)]^{\rm T}$, and $\calK_{:,1}$ denotes the first
column of matrix $\calK$. This is a major advantage of using lattice
points for the construction of the kernel interpolant.

Important properties regarding the kernel interpolant were summarized or
proved in \cite{KKKNS22}:
\begin{itemize}
\item The kernel interpolant is the minimal norm interpolant in the
    sense that it has the \emph{smallest $H$-norm} among all
    interpolants using the same function values of $f$, see
    \cite[Theorem~2.1]{KKKNS22}.

\item The kernel interpolant is optimal in the sense that it has the
    \emph{smallest worst case error} (measured in any norm
    $\|\cdot\|_W$ with $H\subset W$) among all linear or non-linear
    algorithms that use the same function values of $f$, see
    \cite[Theorem~2.2]{KKKNS22}. Recall that the worst case error
    measured in $W$-norm of an algorithm $A$ in $H$ is defined by
\[
  e^{\rm wor}(A;H;W):=\sup_{f\in H,\,\|f\|_H\leq 1}\|f-A(f)\|_W.
\]

\item Any (linear or non-linear) algorithm $A_n$ (with $A_n(0) = 0$)
    that uses function values of $f$ only at the lattice points has
    the lower bound
\[
  e^{\rm wor}(A_n;H;L_p) \ge C\,n^{-\alpha/2}, \quad p\in [1,\infty],
\]
with an explicit constant $C>0$, see \cite[Theorem~3.1]{KKKNS22} and
    \cite{BKUV17}. However, it is known that there exist other
    algorithms based on function values that can achieve an $L_2$
    approximation error upper bound of order $n^{-\alpha}$, see
    \cite{KU21,NSU22,DKU23}. Hence, our lattice-based kernel
    interpolant can only get the half-optimal convergence rate for
    $L_2$ approximation error at best.

\item A component-by-component (CBC) construction from
    \cite{CKNS20,CKNS21} can be used to obtain a lattice generating
    vector such that our lattice-based kernel interpolant achieves
\begin{align} \label{eq:wce}
  e^{\rm wor}(A_n^*;H;L_2) \le \frac{\kappa}{n^{1/(4\lambda)}} \bigg(
	\sum_{\setu\subseteq\{1:s\}}
 \max(|\setu|,1)\,\gamma_{\setu}^\lambda\, [2\zeta(2\alpha\lambda)]^{|\setu|}\bigg)^{1/{(2\lambda)}}
\end{align}
for all $\lambda\in (1/(2\alpha),1]$, with $\kappa :=
\sqrt{2}\,(2.5+2^{4\alpha\lambda+1})^{{{1}/{(4\lambda)}}}$ and
$\zeta(x):=\sum_{k=1}^\infty k^{-x}$ denoting the Riemann zeta
function for $x>1$. Hence, by taking $\lambda\to 1/(2\alpha)$ we
conclude that
\[
  e^{\rm wor}(A_n^*;H;L_2) = \calO(n^{-\alpha/2+\delta}), \quad \delta>0,
\]
where the implied constant depends on $\delta$ but is independent of
    $s$ if the weights $\gamma_\setu$ are such that the sum in
    \eqref{eq:wce} can be bounded independently of $s$, see
    \cite[Theorem~3.2]{KKKNS22}. Consequently, for any $f\in H$ we
    have
\[
  \|f - f_n\|_{L_2} \le e^{\rm wor}(A_n^*;H;L_2)\,\|f\|_H
  = \calO(n^{-\alpha/2+\delta}).
\]
\end{itemize}

The bound \eqref{eq:wce} was proved initially in \cite{CKNS20} only for
prime~$n$, but has since been generalized to composite~$n$ and extensible
lattice sequences in~\cite{KMN}.

Although the theoretical error bound \eqref{eq:wce} holds for all general
weights $\gamma_\setu$, practical implementation of the CBC construction
must take into account the specific form of weights for computational
efficiency. Fast (FFT-based) CBC constructions were developed in
\cite{CKNS21} for product weights, POD weights and SPOD weights with
varying computational cost. Evaluating the kernel interpolant
\eqref{eq:appro} also requires evaluations of the kernel \eqref{eq:ker}
with varying computational cost depending on the form of weights, see
\cite[Section~5.2]{KKKNS22}. Furthermore, if we are interested in
evaluating the kernel interpolant $f_n(\bsy_\ell)$ at $L$ arbitrary points
$\bsy_\ell$, $\ell=1,\ldots,L$, then due to the matrix \eqref{eq:matrix}
being circulant, we can evaluate the kernel interpolant at all the shifted
points $f_n(\bsy_\ell + \bst_{k'})$, $k'=1,\ldots,n$, with only an extra
logarithmic factor in the cost, see \cite[Section~5.1]{KKKNS22}. We
summarize these cost considerations in Table~\ref{tab:cost} (taken from
\cite[Table~1]{KKKNS22}). Clearly, product weights are the most efficient
form of weights in all considerations.

\begin{table}[h]
 \begin{center}
 {
 \ifdefined\journalstyle \else \footnotesize \fi
 \caption{
 \ifdefined\journalstyle \else \footnotesize \fi
 Cost breakdown for the kernel interpolant $f_n$ based on $n$ lattice points $\bst_k$ in $s$ dimensions, evaluated at $L$ arbitrary points $\bsy_\ell$. Here $X$ is the cost for one evaluation
 of $f$.} 
	\label{tab:cost} \setlength{\arrayrulewidth}{0.2pt}
	\medskip\centering
	\begin{tabular}{|l|l@{\,}|l@{\,}|l@{\,}|}
		\hline
		\hfill Operation $\backslash$ Weights & Product & POD & SPOD \\
		\hline
		Fast CBC construction for $\bsz_{\rm gen}$
		& $s\,n\log(n)$
        & $s\,n\log(n) + s^2\log(s)\, n$
		& $s\,n\log(n) + s^3\alpha^2\, n$ \\
		Compute $K(\bst_k,\bszero)$ for all $k$
		& $s\,n$ & $s^2\, n$ & $s^2\,\alpha^2\,n$ \\
		Evaluate $f(\bst_k)$ for all $k$
		& $X\,n$ & $X\,n$ & $X\,n$ \\
		Linear solve for all coefficients $a_k$
		& $n\log(n)$ & $n\log(n)$ & $n\log(n)$ \\
		\hline
		Compute $K(\bst_k,\bsy_\ell)$ for all $k,\ell$
		& $s\,n\,L$ & $s^2\, n\,L$ & $s^2\,\alpha^2\,n\,L$ \\
		Assemble $f_n(\bsy_\ell)$ for all $\ell$
		& $n\,L$ & $n\,L$ & $n\,L$ \\
		OR Assemble $f_n(\bsy_\ell+\bst_{k})$ for all $\ell,k$
		& $n\,\log(n)\,L$ & $n\,\log(n)\,L$ & $n\,\log(n)\,L$ \\
		\hline
	\end{tabular}}
 \end{center}
\end{table}

\section{Kernel interpolant for parametric elliptic PDEs}\label{sec:pde}

In the literature of ``tailored'' quasi-Monte Carlo (QMC) rules for
parametric PDEs, it is customary to analyze the parametric regularity of
the PDE solutions. This information can be used to construct QMC rules
satisfying rigorous error bounds. Many of these studies have been carried
out under the assumption of the so-called ``affine and uniform setting''
as in \eqref{eq:affine}. Examples include the source problem for elliptic
PDEs with random coefficients~\cite{KSS12,DKLNS14,KN16,KN18,GHS18},
spectral eigenvalue problems under
uncertainty~\cite{GGKSS19,GilSch23a,GilSch23b}, Bayesian inverse
problems~\cite{GanPet18,DGLS19,HKS21}, domain uncertainty
quantification~\cite{HPS16}, PDE-constrained optimization under
uncertainty~\cite{GKKSS21,GKKSS}, and many others.
When the input random field is modified to involve a composition with a
periodic function as in \eqref{eq:periodic}, the regularity bound
naturally changes, as we have encountered in \cite{KKS20,KKKNS22,HHKKS}.

We return now to the PDE problem \eqref{eq:notilde_pde} together with the
periodic random field~\eqref{eq:periodic}. We remark that in many studies
the input random field is modeled as an infinite series and the effect of
dimension truncation is analyzed. We will take this point of view below,
as we did in \cite{KKKNS22}. Thus we now have a countably infinite
parameter sequence $\bsy\in U := [0,1]^\bbN$, and $U_s$ in \eqref{eq:u},
\eqref{eq:notilde_pde} and \eqref{eq:periodic} is now replaced by $U$. We
will abuse the notation from the introduction and instead use
$a(\cdot,\bsy)$ and $u(\cdot,\bsy)$ to denote the corresponding random
field and PDE solution, while we write $a_s(\cdot,\bsy) :=
a(\cdot,(y_1,\ldots,y_s,0,0,\ldots))$ and $u_s(\cdot,\bsy) :=
u(\cdot,(y_1,\ldots,y_s,0,0,\ldots,))$ to denote the dimension truncated
counterparts.

Since we have two sets of variables $\bsx\in D$ and $\bsy\in U$, from now
on we will make the domains $D$ and $U$ explicit in our notation.
Let $H_0^1(D)$ denote the subspace of $H^1(D)$ with zero trace on
$\partial D$.  We equip $H_0^1(D)$ with the norm
$\|v\|_{H_0^1(D)}:=\|\nabla v\|_{L_2(D)}$. Let $H^{-1}(D)$ denote the
topological dual space of $H_0^1(D)$, and let $\langle
\cdot,\cdot\rangle_{H^{-1}(D),H_0^1(D)}$ denote the duality pairing
between $H^{-1}(D)$ and $H_0^1(D)$.
We have the parametric weak formulation: for $\bsy\in U$, find
$u(\cdot,\bsy)\in H_0^1(D)$ such that
\begin{align} \label{eq:weak}
 \int_D a(\bsx,\bsy)\nabla u(\bsx,\bsy)\cdot\nabla v(\bsx)\,{\rm d}\bsx
=\langle q,v\rangle_{H^{-1}(D),H_0^1(D)}\quad\text{for all}~v\in H_0^1(D),
\end{align}
where $q\in H^{-1}(D)$. Following the problem formulation in
\cite{KKKNS22}, we make these standing assumptions:
\begin{addmargin}[1.3em]{0em}
\begin{enumerate}
\item[(A1)] $a_0\in L_\infty(D)$ and $\psi_j\in L_\infty(D)$ for all
    $j\geq 1$, and $\sum_{j\ge1} \|\psi_j\|_{L_\infty(D)}<\infty$;
\item[(A2)] there exist $a_{\min}$ and $a_{\max}$ such that
    $0<a_{\min}\leq a(\bsx,\bsy)\leq a_{\max}<\infty$ for all $\bsx\in
    D$ and $\bsy\in U$;
\item[(A3)] $\sum_{j\ge 1} \|\psi_j\|_{L_\infty(D)}^p<\infty$ for some
    $p\in(0,1)$;
\item[(A4)] $a_0\in W^{1,\infty}(D)$ and $\sum_{j\ge 1}
    \|\psi_j\|_{W^{1,\infty}(D)}<\infty$, where
    \ifdefined\journalstyle we defined \fi          
    $\|v\|_{W^{1,\infty}(D)}:=\max\{\|v\|_{L_\infty(D)},\|\nabla
v\|_{L_\infty(D)}\}; $
\item[(A5)] $\|\psi_1\|_{L_\infty(D)}\geq
    \|\psi_2\|_{L_\infty(D)}\geq\cdots$;
\item[(A6)] the spatial domain $D\subset\mathbb R^d$, $d\in\{1,2,3\}$,
    is a convex and bounded polyhedron.
\end{enumerate}
\end{addmargin}

In practical computations, one typically needs to discretize the PDE
\eqref{eq:notilde_pde} using, e.g., a finite element method. While the
weak solution of the PDE problem is in general a Sobolev function and may
not be pointwise well-defined with respect to the spatial variable
$\bsx\in D$, the finite element solution is pointwise well-defined with
respect to $\bsx\in D$ which we may exploit in the construction of the
kernel interpolant. To this end, let $\{V_h\}_h$ be a family of conforming
finite element subspaces $V_h\subset H_0^1(D)$, indexed by the mesh size
$h>0$ and spanned by continuous, piecewise linear finite element basis
functions. Furthermore, we assume that the triangulation corresponding to
each $V_h$ is obtained from an initial, regular triangulation of $D$ by
recursive, uniform partition of simplices. For $\bsy\in U$, the finite
element solution $u_h(\cdot,\bsy)\in V_h$ satisfies
$$
 \int_D a(\bsx,\bsy)\nabla u_h(\bsx,\bsy)\cdot\nabla v_h(\bsx)\,{\rm d}\bsx
=\langle q,v_h\rangle_{H^{-1}(D),H_0^1(D)}\quad\text{for all}~v_h\in V_h.
$$

Let $\bsnu\in\mathbb N_0^\infty$ denote a multi-index with finite order
$|\bsnu| := \sum_{j\ge 1} \nu_j < \infty$, and let $\partial^\bsnu :=
\prod_{j\ge 1} (\partial/\partial y_j)^{\nu_j}$ denote the mixed partial
derivatives with respect to $\bsy$. The standing assumptions (A1)--(A6)
together with the Lax--Milgram lemma ensure that the weak formulation
\eqref{eq:weak} has a unique solution such that for all $\bsy\in U$ (see
\cite{KKS20} for a proof),
\begin{align} \label{eq:reg}
\|\partial^{\bsnu} u(\cdot,\bsy)\|_{H_0^1(D)}\leq
 \frac{\|q\|_{H^{-1}(D)}}{a_{\min}}\,(2\pi)^{|\bsnu|}
 \sum_{\bsm\leq\bsnu}|\bsm|!\,\prod_{j\ge 1} \big( b_j^{m_j}\,S(\nu_j,m_j)\big),
\end{align}
where (no factor $1/\sqrt{6}$ here compared to \cite{KKS20})
\begin{equation}\label{eq:bjdef}
b_j:=\frac{\|\psi_j\|_{L_\infty(D)}}{a_{\min}}\quad\text{for all}~j\geq 1,
\end{equation}
and $S(\nu, m)$ denotes the \emph{Stirling number of the second kind} for
integers $\nu \ge  m\ge 0$, under the convention $S(\nu, 0) = \delta_{\nu,
0}$.

Note that the parametric regularity bound \eqref{eq:reg} holds when
$u(\cdot,\bsy)$ is replaced by a conforming finite element approximation
$u_h(\cdot,\bsy)$. The same is also true for the dimension truncated
solution $u_s(\cdot,\bsy)$ and its corresponding finite element
approximation $u_{s,h}(\cdot,\bsy)$.

Let $H(U_s) = H$ denote the RKHS of functions with respect to $\bsy\in
U_s$ from Section~\ref{section:kernel}. In the framework of
Section~\ref{section:kernel}, for every $\bsx\in D$, we wish to
approximate the dimension truncated finite element solution $f :=
u_{s,h}(\bsx, \cdot)$ at $\bsx$ as a function of $\bsy$, and we define
\[
  f_n:= u_{s,h,n}(\bsx,\cdot):=A_n^*(u_{s,h}(\bsx,\cdot)) \in H(U_s)
\]
to be the corresponding kernel interpolant. Then we are interested in the
joint $L_2$ error
\[
  \| u - u_{s,h,n} \|_{L_2(D\times U)}
  := \sqrt{ \int_D\int_U \left( u(\bsx,\bsy) - u_{s,h,n}(\bsx,\bsy) \right)^2\rd\bsy\rd\bsx},
\]
where we may interchange the order of integration using Fubini's theorem.
Using the triangle inequality, we split this error into three parts, with
$C>0$ in each case denoting a generic constant independent of $s$, $h$ and
$n$:
\begin{enumerate}
\item The \emph{dimension truncation error} satisfies, see
    \cite[Theorem~4.1]{KKKNS22},
\begin{align} \label{eq:dim}
  \|u - u_s\|_{L_2(D\times U)} \le C\,\|q\|_{H^{-1}(D)}\, s^{-(1/p-1/2)},
  \quad\mbox{with $p$ as in (A3)}.
\end{align}
\item The \emph{finite element error} satisfies, see
    \cite[Theorem~4.3]{KKKNS22},
\[
 \|u_s - u_{s,h}\|_{L_2(D\times U)} \le C\,\|q\|_{H^{-1+t}(D)}\,h^{1+t}
 \quad\mbox{as}\quad h\to 0, \quad t\in [0,1].
\]
\item The \emph{kernel interpolation error} satisfies, see
    \cite[Theorem~4.4]{KKKNS22},
\begin{align} \label{eq:interp-err}
   \|u_{s,h}- u_{s,h,n}\|_{L_2(D\times U)} \le
   C\, \|q\|_{H^{-1}(D)}\, n^{-1/(4\lambda)} \,C_s(\lambda)
\end{align}
for all $\alpha\in\N$ and $\lambda\in(\frac{1}{2 \alpha},1]$, where
\begin{align*}
 [C_s(\lambda)]^{2\lambda}&:=\bigg(\sum_{\setu\subseteq\{1:s\}}\max(|\setu|,1)
 \gamma_{\setu}^\lambda[2\zeta(2\alpha\lambda)]^{|\setu|}\bigg)\\
 &\qquad\times\bigg(\sum_{\setu\subseteq\{1:s\}}\frac{1}{\gamma_{\setu}}\bigg(
 \sum_{\bsm_{\setu} \in\{1:\alpha\}^{|\setu|}} |\bsm_{\setu}|!
 \prod_{j\in\setu}(b_j^{m_j}S(\alpha,m_j))\bigg)^2\bigg)^{\lambda}.
\end{align*}
\end{enumerate}

Specifically, the bound \eqref{eq:interp-err} was obtained by writing
\begin{align*}
 \|u_{s,h}- u_{s,h,n}\|_{L_2(D\times U)}^2
 &= \int_D \|u_{s,h}(\bsx,\cdot) -  A_n^*(u_{s,h}(\bsx,\cdot))\|_{L_2(U_s)}^2 \rd\bsx \\
 &\le [e^{\rm wor}(A_n^*;H(U_s);L_2(U_s))]^2 \int_D \|u_{s,h}(\bsx,\cdot)\|_{H(U_s)}^2\rd\bsx.
\end{align*}
The worst case error can be bounded by \eqref{eq:wce}, while the integral
of the squared $H(U_s)$-norm can be bounded by a sum over $\setu\subseteq
\{1:s\}$ involving $\|\partial^\bsnu u_{s,h}(\cdot,\bsy)\|_{H^1_0(D)}^2$
where each $\nu_j$ is $\alpha$ for $j\in\setu$ and is~$0$ otherwise. The
latter $H^1_0(D)$-norm can be bounded by \eqref{eq:reg}, leading
ultimately to the constant $C_s(\lambda)$ in~\eqref{eq:interp-err}.

The difficulty of the parametric PDE problem is largely determined by the
summability exponent $p$ in (A3). We see in \eqref{eq:dim} that the
smaller $p$ is the faster the dimension truncation error decays.
Naturally, the kernel interpolation error \eqref{eq:interp-err} should be
linked with the summability exponent $p$. In this application, the
parameter~$\alpha$ of the RKHS is actually a free parameter for us to
choose, and so are the weights~$\gamma_\setu$. This is more than just a
theoretical exercise: to be able to implement the kernel interpolant in
practice, we must specify $\alpha$ and the weights $\gamma_\setu$, since
they appear in the formula for the kernel \eqref{eq:ker}. The paper
\cite{KKKNS22} proposed a number of choices, all with the aim of
optimizing the convergence rate while keeping the constant $C_s(\lambda)$
bounded independently of $s$. The best convergence rate obtained in
\cite{KKKNS22} was
\begin{align} \label{eq:rate}
  \|u_{s,h} - u_{s,h,n}\|_{L_2(D\times U)} \le C \|q\|_{H^{-1}(D)}\, n^{-r},
  \quad\mbox{with}\quad r = \frac{1}{2p} - \frac{1}{4},
\end{align}
and this was achieved by a choice of SPOD weights. More precisely:
\begin{itemize}
\item We can choose weights (of a very complicated form) to minimize
    $C_s(\lambda)$ and they achieve \eqref{eq:rate}, see
    \cite[Theorem~4.5]{KKKNS22}.\newpage
\item We can choose SPOD weights to mimic the previous weights and
    they also achieve \eqref{eq:rate}, see
    \cite[Theorem~4.5]{KKKNS22}. These SPOD weights are given
    explicitly by
\begin{align} \label{eq:spod}
 \gamma_{\setu}:=\sum_{\bsm_{\setu}\in\{1:\alpha\}^{|\setu|}}(|\bsm_{\setu}|!)^{\frac{2}{1+\lambda}}
 \prod_{j\in\setu}\bigg(\frac{b_j^{m_j}S(\alpha,m_j)}{\sqrt{2{\rm e}^{1/{\rm e}}\zeta(2\alpha\lambda)}}\bigg)^{\frac{2}{1+\lambda}},\\[-.605cm]\notag
\end{align}
with
\begin{align} \label{eq:alpha}
 \alpha:= \left\lfloor \frac{1}{p}+\frac{1}{2}\right\rfloor
 \quad\mbox{and}\quad \lambda:=\frac{p}{2-p}.\\[-.605cm]\notag
\end{align}

\item If $p\in \bigcup_{k=1}^\infty (\frac{2}{2k+1},\frac{1}{k})$ in
    (A3) then we can choose POD weights to achieve \eqref{eq:rate},
    see \cite[Theorem~4.6]{KKKNS22}.

\item We can choose product weights to achieve \eqref{eq:rate} with a
    reduced rate $r$ around one half of $\frac{1}{2p} - \frac{1}{4}$,
    see \cite[Theorem~4.7]{KKKNS22}.
\end{itemize}

Hence, SPOD weights and POD weights can achieve theoretically better
convergence rates than product weights, but they are much more costly as
seen in Table~\ref{tab:cost}. The paper \cite{KKKNS22} reported
comprehensive numerical experiments with the different choices of weights
for the PDE problems of varying difficulties, and concluded that indeed
SPOD weights perform mostly better than POD weights and product weights.
However, the greater computational cost of SPOD weights is definitely
real. We therefore set out to seek better weights.

\section{Seeking better weights}\label{sec:serendipity}

In implementations of the lattice-based kernel interpolant of
\cite{KKKNS22} the weights $\gamma_\setu$ play a dual role.  On the one
hand they appear in the component-by-component (CBC) construction for
determining the lattice generating vector $\bsz_{{\rm gen}}$, which in
turn determines the interpolation points through \eqref{eq:latticepoints}.
On the other hand they appear in the formula \eqref{eq:ker} for the
kernel. In both roles only special forms of weights are feasible, given
that there are $2^s$ subsets $\setu\subseteq \{1:s\}$. The SPOD weights
described above are feasible (and were used in the computations in
\cite{KKKNS22}), but encounter two computational bottlenecks:
\begin{enumerate}
\item The CBC construction used to obtain the generating vector
    $\bsz_{{\rm gen}}$ using SPOD weights has cost $\mathcal
    O(s\,n\log(n) + s^3\,\alpha^2\,n)$, see Table~\ref{tab:cost}
    and~\cite{CKNS20,CKNS21}.
\item Evaluating the kernel interpolant at $L$ arbitrary points using
    SPOD weights has cost $\calO(s^2\,\alpha^2\,n\,L)$, see
    Table~\ref{tab:cost} and \cite[Section~5.2]{KKKNS22}.
\end{enumerate}
While the cost of obtaining the generating vector could be regarded as an
\emph{offline step} that only needs to be performed once for a given set
of problem parameters, the cost of evaluating the kernel interpolant makes
its \emph{online} use unattractive for high-dimensional problems.

We propose the following new formula for product weights to be used in
both roles for the construction of the kernel interpolant:\footnote{Equation~\eqref{eq:serendipity} in the published version of this paper contains a misplaced parenthesis in the final expression for the serendipitous weights. The corrected expression is given here. The numerical results in the paper were
computed with the correct weights and are therefore unaffected. We thank Academia Co-Scientist (Academia.edu) for pointing out this error.}
\begin{align} \label{eq:serendipity}
 \gamma_{\setu}
 :=\sum_{\bsm_{\setu}\in\{1:\alpha\}^{|\setu|}}\prod_{j\in\setu}\bigg(\frac{b_j^{m_j}S(\alpha,m_j)}{\sqrt{2{\rm e}^{1/{\rm e}}\zeta(2\alpha\lambda)}}\bigg)^{\frac{2}{1+\lambda}}
 =\prod_{j\in\setu}\sum_{m=1}^{\alpha}\bigg(\frac{b_j^m S(\alpha,m)}{\sqrt{2{\rm e}^{1/{\rm
 e}}\zeta(2\alpha\lambda)}}\bigg)^{\frac{2}{1+\lambda}},
\end{align}
where $\alpha$ and $\lambda$ are given by \eqref{eq:alpha}. The
weights~\eqref{eq:serendipity} have been obtained from the SPOD
weights~\eqref{eq:spod} simply by leaving out the factorial factor
$(|\bsm_\setu|!)^{2/(1+\lambda)}$. We call these \emph{serendipitous
weights}.

The performance of these weights will be compared against the SPOD
weights~\eqref{eq:spod} in a series of numerical experiments in
Section~\ref{sec:numex}. In addition to the smaller observed errors, the
serendipitous weights (because they are product weights) have obvious
computational advantages:
\begin{enumerate}
\item The CBC construction used to obtain the generating vector
    $\bsz_{\rm gen}$ using product weights has cost $\mathcal
    O(s\,n\log(n))$, see Table~\ref{tab:cost}
    and~\cite{CKNS20,CKNS21}.
\item Evaluating the kernel interpolant at $L$ arbitrary points using
    product weights has cost $\mathcal O(s\,n\,L)$, see
    Table~\ref{tab:cost} and~\cite[Section~5.2]{KKKNS22}.
\end{enumerate}
As we shall see, serendipitous weights~\eqref{eq:serendipity} allow us to
tackle successfully very high-dimensional approximation problems.
Moreover, we still have the rigorous error bound given in
\eqref{eq:interp-err}. We no longer have a guarantee of a constant
independent of $s$, but the observed performance will be seen to be
excellent.

\section{Numerical experiments}\label{sec:numex}

We consider the parametric PDE problem~\eqref{eq:pde} converted to
periodic form in \eqref{eq:notilde_pde}, with the periodic diffusion
coefficient \eqref{eq:periodic}. The domain is $D=(0,1)^2$, and the source
term is $q(\bsx)=x_2$. For the mean field we set $a_0(\bsx) = 1$, and for
the stochastic fluctuations we take the functions
\[ 
\psi_j(\bsx):=cj^{-\theta}\sin(j\pi x_1)\sin(j\pi x_2),\quad \bsx=(x_1,x_2)\in D,~j\geq 1,
\] 
where $c>0$ is a constant determining the magnitude of the fluctuations,
and $\theta>1$ is the decay rate of the stochastic fluctuations. The
sequence $(b_j)_{j\geq1}$ defined by \eqref{eq:bjdef} becomes
\begin{align} \label{eq:bj}
  b_j:=\frac{cj^{-\theta}}{a_{\min}}, \quad j\ge 1,
\end{align}
where for simplicity we take
\[
  a_{\min}:=1-c\zeta(\theta) \quad\mbox{and}\quad
  a_{\max}:=1+c\zeta(\theta),
\]
and enforce $c<\frac{1}{\zeta(\theta)}$ to ensure the uniform ellipticity
condition.

For each fixed $\bsy \in [0,1]^s$ we solve the PDE using a piecewise
linear finite element method with $h = 2^{-5}$ as the finite element mesh
size. We construct a kernel interpolant
$u_{s,h,n}(\bsx,\cdot)=A_n^*(u_{s,h}(\bsx,\cdot))$ for the finite element
solution $u_{s,h}(\bsx,\cdot)$ using both the SPOD weights~\eqref{eq:spod}
and serendipitous weights~\eqref{eq:serendipity}. These weights also enter
the CBC construction used to obtain a lattice generating vector satisfying
\eqref{eq:interp-err}: specifically, the kernel interpolant is constructed
over the point set $\bst_k=\{k\bsz_{{\rm gen}}/n\}$, where $\bsz_{{\rm
gen}}$ is obtained using the algorithm described in~\cite{CKNS21}.

The kernel interpolation error is estimated by computing
\begin{align}\label{eq:errorbound}
{\rm error}&=\sqrt{\int_D\int_{[0,1]^s}
 \big(u_{s,h}(\bsx,\bsy)-u_{s,h,n}(\bsx,\bsy)\big)^2\rd\bsy\rd\bsx}\nonumber\\
&\approx\sqrt{\frac{1}{Ln}\sum_{\ell=1}^L \sum_{k=1}^n
\int_D \big(u_{s,h}(\bsx,\bsy_{\ell}+\bst_k)-u_{s,h,n}(\bsx,\bsy_\ell+\bst_k)\big)^2\,{\rm d}\bsx},
\end{align}
where $\bsy_{\ell}$ for $\ell=1,\ldots,L$ denotes the sequence of Sobol$'$
points with $L=100$. Note that since the functions $u_{s,h}(\bsx,\bsy)$
and $u_{s,h,n}(\bsx,\bsy)$ are 1-periodic with respect to~$\bsy$, the
kernel interpolant can be evaluated efficiently over the union of shifted
points $\bsy_{\ell}+\bst_k$ for $\ell=1,\ldots,L$ and $k=1,\ldots,n$,
using FFT, see Table~\ref{tab:cost} and~\cite[Section~5.1]{KKKNS22}.

\subsection{Fixing the parameters in the weights}

To implement the kernel interpolant with either SPOD or serendipitous
weights, one first has to choose the parameters $c$ and $\theta$ in
\eqref{eq:bj}. The next step is to decide on a value of $p\in (0,1)$ that
satisfies (A3). This clearly restricts $p$ to the smaller interval $p\in
(1/\theta, 1).$  The choice of $p$ in turn determines the parameters
$\alpha$ and $\lambda$ through~\eqref{eq:alpha}.

In the experiments we choose three different values for the decay
parameter~$\theta$, namely, $\theta = 3.6$, $2.4$ and~$1.2$, ranging from
the very easy to the very difficult. Correspondingly, we choose $p =
\frac{1}{3.3}$, $\frac{1}{2.2}$ and~$\frac{1}{1.1}$, respectively, leading
to values of the smoothness parameter $\alpha = 3$, $2$ and $1$,
respectively. We also use different values for the parameter $c =
\frac{0.2}{\sqrt 6}$, $\frac{0.4}{\sqrt 6}$, and $\frac{1.5}{\sqrt 6}$,
again ranging from easy to difficult. (The factor $\frac{1}{\sqrt 6}$ has
been included here to facilitate comparisons to the numerical results
in~\cite{KKKNS22}.)

\subsection{Comparing SPOD weights with serendipitous weights}

\begin{figure}[!t]
 \ifdefined\journalstyle
 \includegraphics[height=.5\textwidth]{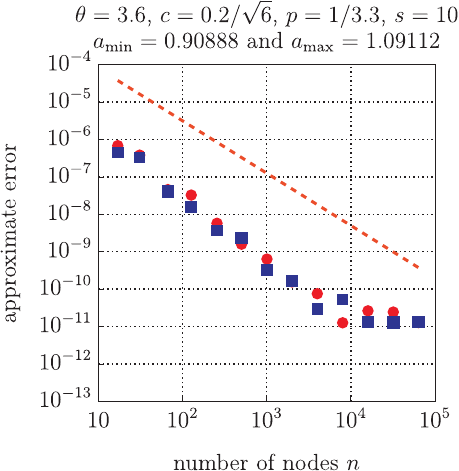}
 \includegraphics[height=.5\textwidth,clip,trim=0cm 0cm 6.3cm 0cm]{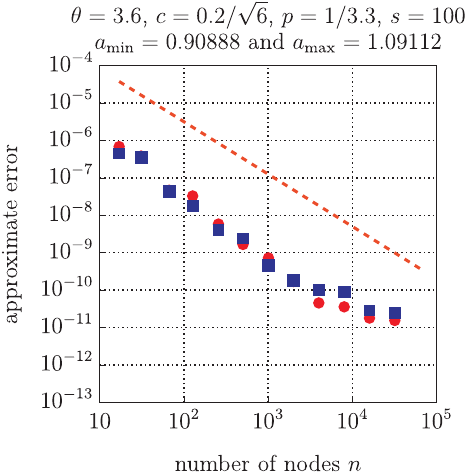}\\
 \\
 \includegraphics[height=.5\textwidth]{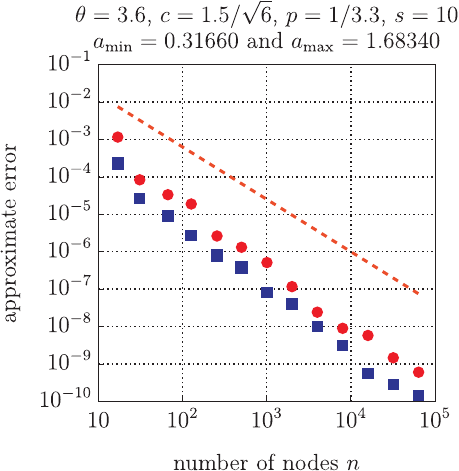}
 \includegraphics[height=.52\textwidth,clip,trim=0cm 0cm 6.3cm 0cm]{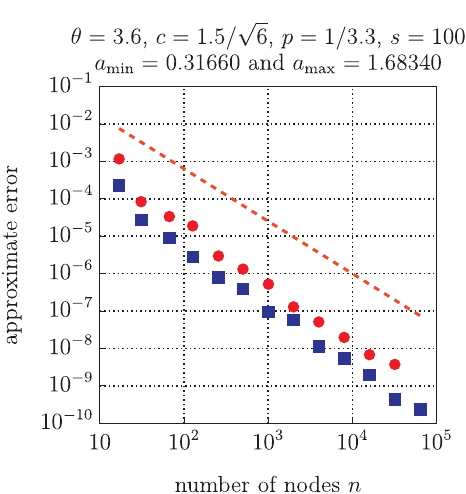}\\
 \centering
 \includegraphics[width=.32\textwidth,clip,trim=8.9cm 4.59cm 0cm 2.7cm]{figures/err_s100_6X.pdf}
 \includegraphics[width=.32\textwidth,clip,trim=8.9cm 3.75cm 0cm 3.5cm]{figures/err_s100_6X.pdf}
 \includegraphics[width=.32\textwidth,clip,trim=8.9cm 2.74cm 0cm 4.6cm]{figures/err_s100_6X.pdf}
 \else
  \includegraphics[height=.4\textwidth]{figures/err_s10_5X.pdf}
 \includegraphics[height=.4\textwidth,clip,trim=0cm 0cm 6.3cm 0cm]{figures/err_s100_5X.pdf}\\
 \\
 \includegraphics[height=.4\textwidth]{figures/err_s10_6X.pdf}
 \includegraphics[height=.42\textwidth,clip,trim=0cm 0cm 6.3cm 0cm]{figures/err_s100_6X.pdf}\\
 \centering
 \includegraphics[width=.28\textwidth,clip,trim=8.9cm 4.59cm 0cm 2.7cm]{figures/err_s100_6X.pdf}
 \includegraphics[width=.28\textwidth,clip,trim=8.9cm 3.75cm 0cm 3.5cm]{figures/err_s100_6X.pdf}
 \includegraphics[width=.28\textwidth,clip,trim=8.9cm 2.74cm 0cm 4.6cm]{figures/err_s100_6X.pdf}
 \fi
 \caption{
 \ifdefined\journalstyle \else \footnotesize \fi
 The kernel interpolation errors of the PDE problem~\eqref{eq:notilde_pde}
 and~\eqref{eq:periodic} with $\theta=3.6$, $p=1/3.3$,
 $c\in\{\frac{0.2}{\sqrt 6},\frac{1.5}{\sqrt 6}\}$, and $s\in\{10,100\}$.
 The results are displayed for kernel interpolants constructed using SPOD
 weights~\eqref{eq:spod} and serendipitous
 weights~\eqref{eq:serendipity}.} \label{fig:easy}
\end{figure}

\begin{figure}[!t]
 \ifdefined\journalstyle
 \includegraphics[height=.5\textwidth]{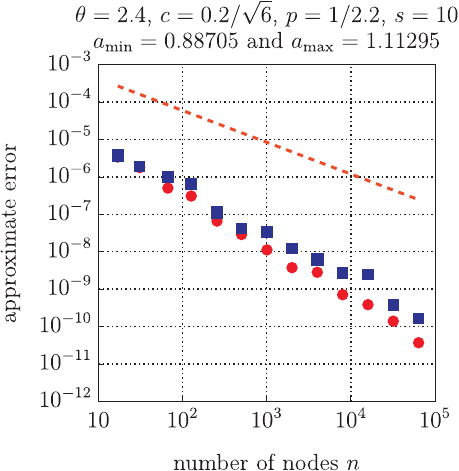}\!\!\!
 \includegraphics[height=.5\textwidth,clip,trim=0cm 0cm 6.3cm 0cm]{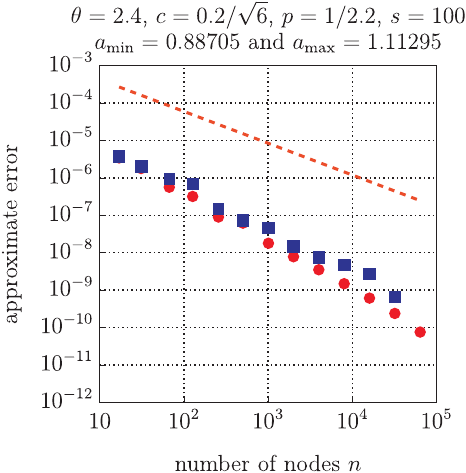}\\
 \\ \vspace*{.1cm}
 \includegraphics[height=.5\textwidth]{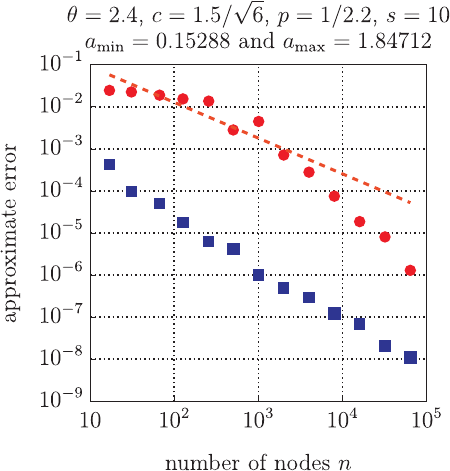}
 \includegraphics[height=.52\textwidth,clip,trim=0cm 0cm 6.3cm 0cm]{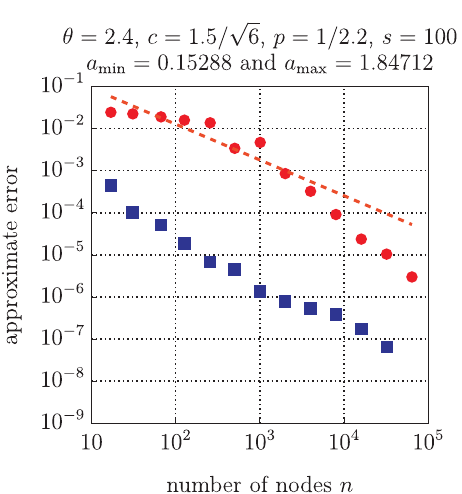}\\
 \centering
 \includegraphics[width=.32\textwidth,clip,trim=8.9cm 4.59cm 0cm 2.7cm]{figures/err_s100_4X.pdf}
 \includegraphics[width=.32\textwidth,clip,trim=8.9cm 3.75cm 0cm 3.5cm]{figures/err_s100_4X.pdf}
 \includegraphics[width=.32\textwidth,clip,trim=8.9cm 2.74cm 0cm 4.6cm]{figures/err_s100_4X.pdf}
 \else
 \includegraphics[height=.4\textwidth]{figures/err_s10_3X.pdf}\!\!\!
 \includegraphics[height=.4\textwidth,clip,trim=0cm 0cm 6.3cm 0cm]{figures/err_s100_3X.pdf}\\
 \\ \vspace*{.1cm}
 \includegraphics[height=.4\textwidth]{figures/err_s10_4X.pdf}
 \includegraphics[height=.42\textwidth,clip,trim=0cm 0cm 6.3cm 0cm]{figures/err_s100_4X.pdf}\\
 \centering
 \includegraphics[width=.28\textwidth,clip,trim=8.9cm 4.59cm 0cm 2.7cm]{figures/err_s100_4X.pdf}
 \includegraphics[width=.28\textwidth,clip,trim=8.9cm 3.75cm 0cm 3.5cm]{figures/err_s100_4X.pdf}
 \includegraphics[width=.28\textwidth,clip,trim=8.9cm 2.74cm 0cm 4.6cm]{figures/err_s100_4X.pdf}
 \fi
 \caption{
 \ifdefined\journalstyle \else \footnotesize \fi
 The kernel interpolation errors of the PDE problem~\eqref{eq:notilde_pde} and~\eqref{eq:periodic}
 with $\theta=2.4$, $p=1/2.2$, $c\in\{\frac{0.2}{\sqrt 6},\frac{1.5}{\sqrt 6}\}$, and $s\in\{10,100\}$.
 The results are displayed for kernel interpolants constructed using SPOD weights~\eqref{eq:spod}
 and serendipitous weights~\eqref{eq:serendipity}.} \label{fig:mid}
\end{figure}

\begin{figure}[!t]
 \ifdefined\journalstyle
 \includegraphics[height=.5\textwidth]{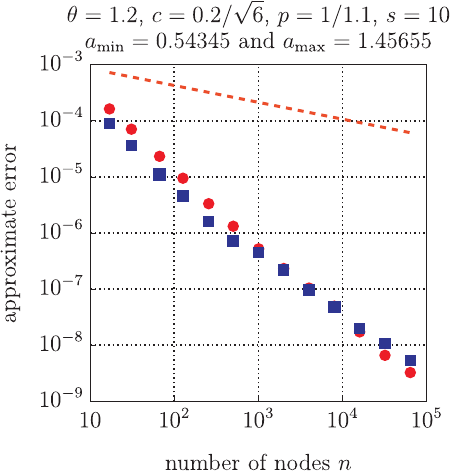}
 \includegraphics[height=.5\textwidth,clip,trim=0cm 0cm 6.3cm 0cm]{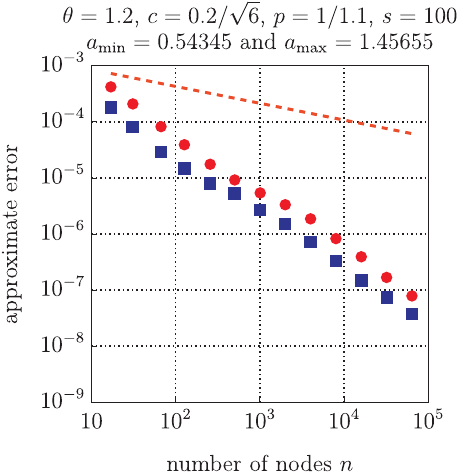}\\
 \\
 \includegraphics[height=.5\textwidth]{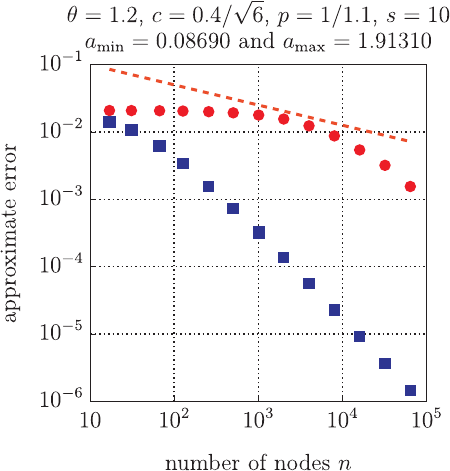}
 \includegraphics[height=.52\textwidth,clip,trim=0cm 0cm 6.3cm 0cm]{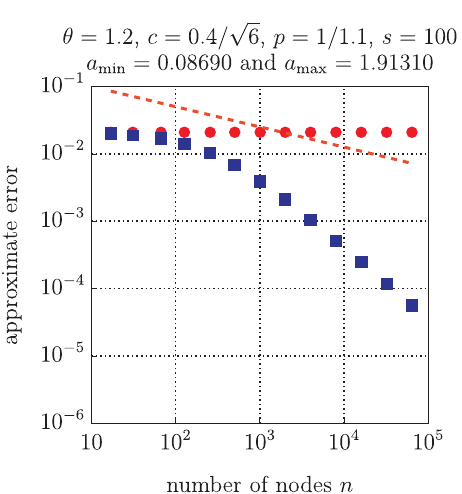}\\
 \centering
 \includegraphics[width=.32\textwidth,clip,trim=8.9cm 4.59cm 0cm 2.7cm]{figures/err_s100_2X.pdf}
 \includegraphics[width=.32\textwidth,clip,trim=8.9cm 3.75cm 0cm 3.5cm]{figures/err_s100_2X.pdf}
 \includegraphics[width=.32\textwidth,clip,trim=8.9cm 2.74cm 0cm 4.6cm]{figures/err_s100_2X.pdf}
 \else
  \includegraphics[height=.4\textwidth]{figures/err_s10_1X.pdf}
 \includegraphics[height=.4\textwidth,clip,trim=0cm 0cm 6.3cm 0cm]{figures/err_s100_1X.pdf}\\
 \\
 \includegraphics[height=.4\textwidth]{figures/err_s10_2X.pdf}
 \includegraphics[height=.42\textwidth,clip,trim=0cm 0cm 6.3cm 0cm]{figures/err_s100_2X.pdf}\\
 \centering
 \includegraphics[width=.28\textwidth,clip,trim=8.9cm 4.59cm 0cm 2.7cm]{figures/err_s100_2X.pdf}
 \includegraphics[width=.28\textwidth,clip,trim=8.9cm 3.75cm 0cm 3.5cm]{figures/err_s100_2X.pdf}
 \includegraphics[width=.28\textwidth,clip,trim=8.9cm 2.74cm 0cm 4.6cm]{figures/err_s100_2X.pdf}
 \fi
 \caption{
 \ifdefined\journalstyle \else \footnotesize \fi
 The kernel interpolation errors of the PDE problem~\eqref{eq:notilde_pde} and~\eqref{eq:periodic}
 with $\theta=1.2$, $p=1/1.1$, $c\in\{\frac{0.2}{\sqrt 6},\frac{0.4}{\sqrt 6}\}$, and $s\in\{10,100\}$.
 The results are displayed for kernel interpolants constructed using SPOD weights~\eqref{eq:spod}
 and serendipitous weights~\eqref{eq:serendipity}.} \label{fig:hard}
\end{figure}

We here compare the kernel interpolation errors using both the SPOD
weights~\eqref{eq:spod} and serendipitous weights~\eqref{eq:serendipity}.
The computed quantity in each case is the estimated $L_2$ error with
respect to both the domain and stochastic variables, given by
\eqref{eq:errorbound}.

The results are displayed in Figures~\ref{fig:easy}, \ref{fig:mid} and
\ref{fig:hard} for the three different $\theta$ values, ranging from the
easiest to the hardest. In each case the graphs on the left are for
dimension $s = 10$, those on the right for $s = 100$.  Each figure also
gives the computed errors for two different values of the parameter $c$,
with the easier (i.e., smaller) value used in the upper pair, the harder
(i.e., larger) value in the lower pair. Each graph also shows (dashed
line) the theoretical convergence rate \eqref{eq:rate} for the given value
of $p$.

The striking fact is that  the serendipitous weights perform about as well
as SPOD weights for all the easier cases (all graphs in
Figure~\ref{fig:easy}, the upper graphs in
Figures~\ref{fig:mid}~and~\ref{fig:hard}),  while dramatically
outperforming the SPOD weights for the harder cases (the lower graphs in
Figures~\ref{fig:mid} and~\ref{fig:hard}).

One way to assess the hardness of a particular parameter set is to inspect
the values of $a_{{\rm min}}$ and $a_{{\rm max}}$ given in the legend
above each graph. In particular, the hardest problem is the fourth graph
in Figure~\ref{fig:hard}, where the dimensionality is $s = 100$, and the
random field has values ranging from around 0.1 to near~2.  For this case
the SPOD weights are seen to fail completely. Yet even in this case the
serendipitous weights perform superbly.

The plateau in the convergence graph for the SPOD weights in
Figure~\ref{fig:hard} can be explained as follows: the SPOD weights in
this case become very large with increasing dimension and, in consequence,
the kernel interpolant becomes very spiky at the lattice points and near
zero elsewhere. Thus we are effectively seeing just the $L_2$~norm of the
target function $u_{s,h}$ for feasible values of $n$.

The intuition behind the serendipitous weights is that the problem of
overly large weights is overcome by omitting the factorials in the SPOD
weight formula \eqref{eq:spod}.

Finally, it is worth emphasizing that the construction cost with
serendipitous weights is considerably cheaper than with SPOD weights. Yet
the quality of the kernel interpolant appears to be just as good as
or---as seen in Figures~\ref{fig:mid}
and~\ref{fig:hard}---\emph{dramatically better than} the interpolant based
on SPOD weights.  Putting these two aspects together, we resolved to
repeat the experiments for a still larger value of $s$, well beyond the
reach of SPOD weights.

\subsection{1000-dimensional examples} 

\begin{figure}[!h]
 \begin{center}
 \ifdefined\journalstyle
 \includegraphics[height=.45\textwidth]{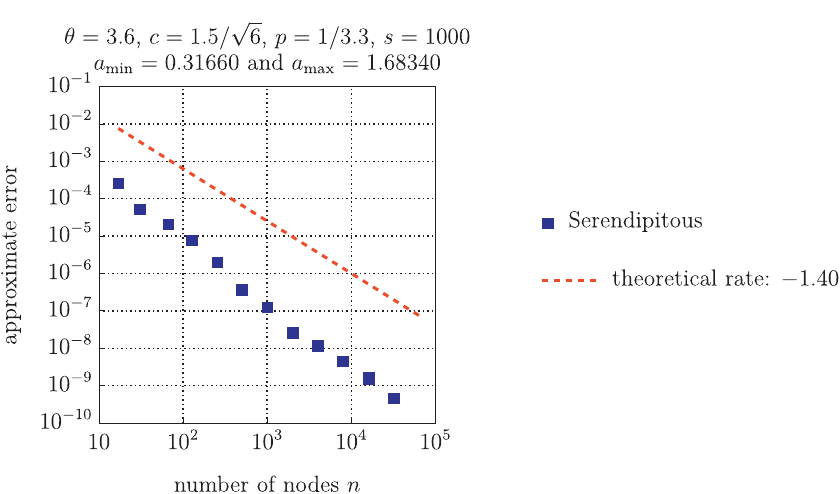}\\
 \includegraphics[height=.45\textwidth]{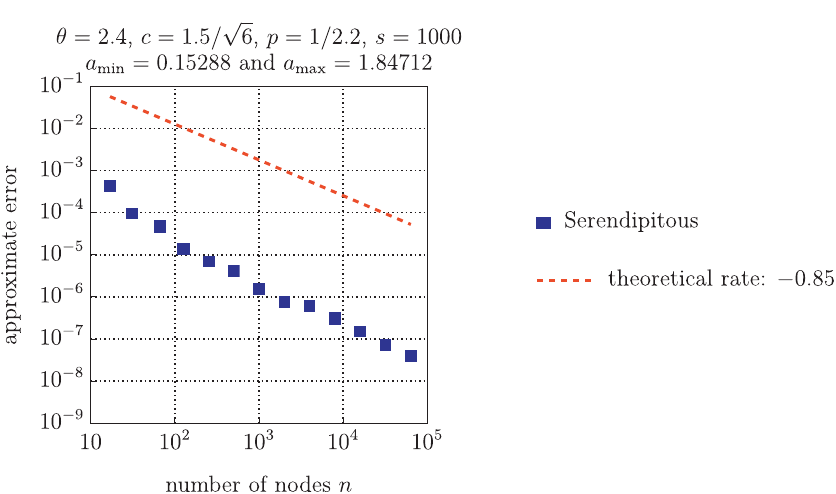}\\
 \includegraphics[height=.45\textwidth]{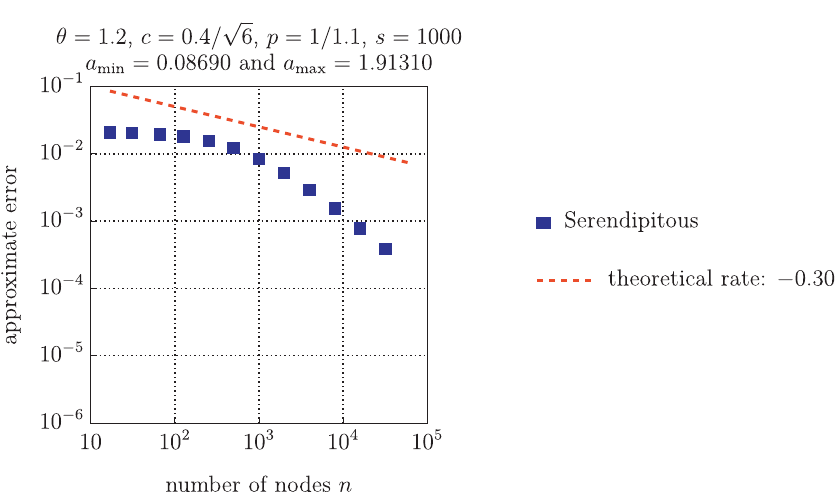}
 \else
 \includegraphics[height=.4\textwidth]{figures/err_s1000_6X.pdf}\\
 \includegraphics[height=.4\textwidth]{figures/err_s1000_4X.pdf}\\
 \includegraphics[height=.4\textwidth]{figures/err_s1000_2X.pdf}
 \fi
 \end{center}
 \caption{
 \ifdefined\journalstyle \else \footnotesize \fi
 The kernel interpolation errors of the PDE problem~\eqref{eq:notilde_pde} and~\eqref{eq:periodic}
 with $s=1000$ and parameters $\theta=3.6$, $c=\frac{1.5}{\sqrt 6}$, and $p=1/3.3$ (top), $\theta=2.4$,
 $c=\frac{1.5}{\sqrt{6}}$, and $p=1/2.2$ (middle), $\theta=1.2$, $c=\frac{0.4}{\sqrt 6}$, and $p=1/1.1$ (bottom).
 The results are displayed for kernel interpolants constructed using serendipitous
 weights~\eqref{eq:serendipity}.} \label{fig:tough}
\end{figure}

Since the serendipitous weights~\eqref{eq:serendipity} are product
weights, we are able to carry out computations using much higher
dimensionalities than before. We consider the previous set of three
$\theta$ values together with the harder value of $c$ in each case, and
now set the upper limit of the series~\eqref{eq:periodic} to be $s=1000$.
The results are displayed in Figure~\ref{fig:tough}.

The computational performance of the kernel interpolant using serendips is
seen in Figure 4 to continue to be excellent even when $s=1000$. The
method works well even for the most difficult experiment, illustrated on
the bottom of Figure~\ref{fig:tough}. While the pre-asymptotic regime is
even longer than in the case $s=100$ (shown in the bottom right of
Figure~\ref{fig:hard}), the kernel interpolation error does not stall as
it does in the case $s=100$ for SPOD weights. Thus the kernel interpolant
based on serendipitous weights appears to be robust in practice.

\ifdefined\journalstyle \else \clearpage \fi

\section{Conclusions}

We have introduced a new class of product weights, called the
\emph{serendipitous weights}, to be used in conjunction with the
lattice-based kernel interpolant presented in~\cite{KKKNS22}. Numerical
experiments illustrate that this family of weights appears to be robust
when it comes to kernel interpolation of parametric PDEs.

Numerical experiments in the paper comparing the performance with
previously studied SPOD weights show that not only are the new weights
cheaper and easier to work with, but also they give much better results
for hard problems.

\section*{Acknowledgements}
F.~Y. Kuo and I.~H. Sloan acknowledge support from the Australian Research
Council (DP210100831). This research includes computations using the
computational cluster Katana~\cite{katana} supported by Research
Technology Services at UNSW Sydney.

\end{document}